\documentclass[12pt]{amsart}
\usepackage{amssymb}
\usepackage{epsfig}
\usepackage{mathabx}
\usepackage{txfonts}
\usepackage{amsmath}
\usepackage{amssymb}
\usepackage{color}
\usepackage{amsmath}
\usepackage{cite}

\numberwithin{equation}{section}

\begin{document}

\newtheorem{mthm}{Theorem}
\newtheorem{mcor}{Corollary}
\newtheorem{mpro}{Proposition}
\newtheorem{mfig}{figure}
\newtheorem{mlem}{Lemma}
\newtheorem{mdef}{Definition}
\newtheorem{mrem}{Remark}
\newtheorem{mpic}{Picture}
\newtheorem{rem}{Remark}[section]
\newcommand{\ra}{{\mbox{$\rightarrow$}}}
\newtheorem{Remark}{Remark}[section]
\newtheorem{thm}{Theorem}[section]
\newtheorem{pro}{Proposition}[section]
\newtheorem*{proA}{Proposition A}
\newtheorem*{proB}{Proposition B}
\newtheorem*{proC}{Proposition C}
\newtheorem*{proD}{Proposition D}
\newtheorem*{proE}{Proposition E}

\newtheorem{lem}{Lemma}[section]
\newtheorem{defi}{Definition}[section]
\newtheorem{cor}{Corollary}[section]

\def \p{\partial}
\def \R{\mathbb{R}}

\title[]{Remark on a special class of Finsler $p$-Laplacian equation}

\author{Yuan Li }
\address{School of Mathematical Sciences,  Key Laboratory of MEA(Ministry of Education) \& Shanghai Key Laboratory of PMMP,  East China Normal University, Shanghai 200241, China}
\email{yli@math.ecnu.edu.cn}

\author{Dong Ye }
\address{ School of Mathematical Sciences,  Key Laboratory of MEA(Ministry of Education) \& Shanghai Key Laboratory of PMMP,  East China Normal University, Shanghai 200241, China and
IECL, UMR 7502, University of Lorraine, 57050 Metz, France}
\email{dye@math.ecnu.edu.cn, dong.ye@univ-lorraine.fr}
\thanks{}

\date{}

\begin{abstract}
We investigate the anisotropic elliptic equation $-\Delta_p^H u = g(u)$. Recently, Esposito, Riey, Sciunzi, and Vuono introduced an anisotropic Kelvin transform in their work \cite{ERSV2022} under the $(H_M)$ condition, where $H(\xi)=\sqrt{\langle M\xi,\xi\rangle}$ with a positive definite symmetric matrix $M$. Here, we emphasize that under the $(H_M)$ assumption, the Finsler $p$-Laplacian and the classical $p$-Laplacian operator are equivalent following a linear transformation. This equivalence offers us a more direct route to derive the pivotal findings presented in \cite{ERSV2022}. While this equivalence is crucial and noteworthy, to our knowledge, it has not been explicitly stated in the current literature.
\end{abstract}

\maketitle


\section{Introduction and main results}
In this note, we consider the following Finsler $p$-Laplacian equation
\begin{align}\label{main eq1}
-\Delta_{p}^{H}u =g(u) \quad \mbox{in }\; \mathbb{R}^{n},
\end{align}
where $p > 1$, $n \geq 2$, $g\in C^{1}(\mathbb{R})$ and the operator $\Delta_{p}^{H}$ is defined by
$$\Delta_{p}^{H}u:=\sum_{i=1}^{n}\frac{\partial}{\partial x_{i}}\left(H^{p-1}(\nabla u)H_{\xi_{i}}(\nabla u)\right),$$
where $H_{\xi_i} = \frac{\p H}{\p \xi_i}$ with $H$, a convex function over $\R^n$ satisfying
\begin{itemize}
\item[-] $H \in C^{2}$ in $\mathbb{R}^{n}\setminus\{0\};$
\item[-] $H(t\xi)=|t|H(\xi)$ for any $t\in\mathbb{R}$, $\xi \in\mathbb{R}^{n};$
\item[-] there exist $0<a\leq b<\infty$ such that $a|\xi|\leq H(\xi)\leq b|\xi|$ for any $\xi \in \mathbb{R}^{n}$.
\end{itemize}
Let $H^{*}$ be the support function of $K:=\{x\in\mathbb{R}^{n}: H(x)<1\}$ defined by
$$H^{*}(x):=\sup_{\xi\in K} \langle x,\xi\rangle.$$
$H^{*}$ is called the dual norm of $H$, and $B_{r}(a):=\{\xi\in\mathbb{R}^{n}: H^{*}(\xi-a)<r\}$ is known as the Wulff ball of radius $r$ centered at $a \in \R^n$.

\smallskip
The anisotropic operators $\Delta_p^H$ have garnered significant attention in the literature. In 1901, Wulff \cite{W1901} pioneered the use of these operators to study crystal shapes and the minimization of anisotropic surface tensions. For a deeper understanding of anisotropic operators and the intricate world of Finsler geometry, we refer to \cite{AFTL1997,BCS2000,CS2006,CS2009,FM1991,SS2016} and references therein.

\smallskip
The most well known example is $$H(\xi)= |\xi|_q := \Big(\sum\limits_{i=1}^{n}|\xi_{i}|^q\Big)^{\frac{1}{q}} \;\; \mbox{ in }\; \R^n$$ for $1<q<\infty$. In particular, when $q=2$, i.e.~$H(\xi)= |\xi|$ is the Euclidean norm, $\Delta_p^H$ is just the classical isotropic $p$-Laplacian $\Delta_p$.

\smallskip
Here we mainly consider two cases: $p=2$ and $p=n$. When $p=2$, we denote $\Delta_{p}^{H}$ as $\Delta^{H}$, which is often called by the Finsler Laplacian or anisotropic Laplacian.
To consider \eqref{main eq1}, let us precise first the meaning of weak solution.
\begin{defi}
\label{def1}
Let $f \in L^{1}_{loc}(\mathbb{R}^{n})$ and $p > 1$. We say that $u$ is a weak solution of $-\Delta_p^H u = f$ in $\R^n$, if $u\in W^{1,p}_{loc}(\mathbb{R}^{n})$ and
$$\int_{\mathbb{R}^{n}}H^{p-1}(\nabla u)H_{\xi}(\nabla u)\nabla\phi dx=\int_{\mathbb{R}^{n}}f\phi dx, \quad \mbox{for all }\; \phi\in C_{c}^1(\mathbb{R}^{n}).$$
In particular, we say that $u$ is a weak solution of (\ref{main eq1}) if $f = g(u)$.
\end{defi}

Ciraolo, Figalli and Roncoroni \cite{CFR2020} observed that generally the Kelvin transform does not work in anisotropic case.
Very recently, Esposito, Riey, Sciunzi and Vuono considered in \cite{ERSV2022} some special $H$ satisfying the so-called $(H_M)$ condition, that is,
\begin{align}\label{hm}
H(\xi)=\sqrt{\langle M\xi,\xi\rangle}, \quad \mbox{for any } \; \xi \in \R^n,
\end{align}
where $M$ is a positive definite symmetric matrix and $n \geq 2$. In this case, it's easy to see that
\begin{align}\label{hm'}
H^{*}(\xi)=\sqrt{\langle M^{-1}\xi,\xi\rangle},
\end{align}
where $M^{-1}$ is the inverse matrix of $M$. In \cite{ERSV2022}, they obtained an anisotropic version of Kelvin transform as follows.
\begin{align}\label{Kelvin}
T_{H}: \mathbb{R}^{n}\setminus\{0\}\rightarrow\mathbb{R}^{n}\setminus\{0\}, \quad T_{H}(\xi):=\frac{\nabla H(\xi)}{H(\xi)},
\end{align}
and
$$\widehat{u}:=\frac{u\circ T_{H}}{H^{n-2}} \quad\mbox{in}\quad\mathbb{R}^{n}\setminus\{0\}.$$
More precisely, they proved (see \cite[Theorems 1.2]{ERSV2022})
\begin{proA}\label{pro1}
Let $H$ satisfy (\ref{hm}). If $u\in H^{1}_{loc}(\mathbb{R}^{n})$ is a locally bounded solution of $-\Delta^{H}u=f(x)$ in $\mathbb{R}^{n}$ with $f\in L^{2}_{loc}(\mathbb{R}^{n})$, then $\widehat{u}\in H^{1}_{loc}(\mathbb{R}^{n}\setminus\{0\})$ and weakly solves the dual equation
$$-\Delta^{H^{*}}\widehat{u}=\frac{f\circ T_{H}}{H^{n+2}} \quad\mbox{in}\quad\mathbb{R}^{n}\setminus\{0\}.$$
\end{proA}

They also studied the Kelvin transform for $-\Delta_n^H$ in $\R^n$ (see \cite[Theorems 1.4]{ERSV2022}) under the $(H_M)$ assumption.
\begin{proB}\label{pro3}
Let $H$ satisfy (\ref{hm}). If $u\in W^{1,n}_{loc}(\mathbb{R}^{n})$ is a locally bounded solution of $-\Delta^{H}_{n}u=g(x)$ in $\mathbb{R}^{n}$ with $g\in L^\frac{n}{n-1}_{loc}(\mathbb{R}^{n})$, then $u^{*}:=u\circ T_{H}$ with $T_H$ given by \eqref{Kelvin} weakly solves the dual equation
$$-\Delta^{H^{*}}_{n}u^{*}=\frac{g\circ T_{H}}{H^{2n}} \quad\mbox{in}\quad\mathbb{R}^{n}\setminus\{0\}.$$
\end{proB}

It is worthy to mention that Ferone and Kawohl \cite{FK2009} first put forward the following assumptions for $H$, that is
\begin{align}\label{hm1}
\langle H(x)\nabla H(x),H^{*}(y)\nabla H^{*}(y)\rangle=\langle x,y\rangle, \quad \forall\, x,y\in\mathbb{R}^{n}.
\end{align}
Assuming \eqref{hm1}, they proved that the mean value property for anisotropic harmonic functions, i.e. solutions of the equation
$$\Delta^{H}u=0.$$
\begin{proC}\label{pro4}
Suppose that $H$ and $H^{*}$ satisfy (\ref{hm1}). If $\Delta^{H}u=0$ in $B_{\rho}(x_{0})$, then for every Wulff ball of radius $r\in(0,\rho)$, $u$ satisfies
$$u(x_{0})=\frac{1}{n\kappa r^{n-1}}\int_{\partial B_{r}(x_{0})}u(x)dS = \frac{1}{\kappa r^{n}}\int_{ B_{r}(x_{0})}u(x)dx.$$
Here $\kappa$ stands for the Lebesgue measure of the unit Wulff ball $B_1(0)$.
\end{proC}

However, a pertinent question arises: which functions $H$ can satisfy the condition (\ref{hm1})? It is immediately evident that the assumption \eqref{hm} inherently fulfills (\ref{hm1}). Interestingly, Cozzi, Farina, and Valdinoci \cite{CFV2016} revealed a profound criterion: When $K =\{x\in\mathbb{R}^{n}: H(x)<1\}$ is strictly convex, (\ref{hm1}) and (\ref{hm}) are indeed equivalent.
\begin{proD}\label{pro8}
Let $H\in C^{1}(\mathbb{R}^{n}\setminus\{0\})$ be a positive homogeneous function of degree $1$ satisfying $H(\xi)>0$ for all $\xi\ne 0$. Assume that $K$ is strictly convex. Then the assumption (\ref{hm1})
is equivalent to asking $H$ to satisfy $(H_M)$, i.e. $H(\xi)=\sqrt{\langle M\xi,\xi\rangle}$, for some symmetric and positive definite matrix $M$.
\end{proD}

Since $M$ is a symmetric and positive definite matrix, it is well-known that there exists a unique symmetric and positive definite matrix $B$ satisfying $B^2=M$, or equivalently $B = \sqrt{M}$.

\medskip
Utilizing the linear transformation $x \mapsto Bx$, we show an elementary fact: under the $(H_M)$ assumption, the operator $\Delta^{H}_{p}$ exhibits a close relationship with the classical $p$-Laplacian operator. Consequently, Propositions A through C can be derived directly from their analogous conclusions in isotropic settings.
\begin{thm}\label{thm1}
Suppose that $H$ satisfies (\ref{hm}). Then for any $p >1$ and any $W^{1, p}_{loc}$ function $u$, if $\widetilde{u}(x)=u(Bx)$,  there holds
\begin{align}
\label{equa1}
\int_{\mathbb{R}^{n}}|\nabla \widetilde{u}|^{p-2}\nabla \widetilde{u}\cdot\nabla\widetilde{\phi}dx=\int_{\mathbb{R}^{n}}\left[H^{p-1}(\nabla u)H_{\xi}(\nabla u)\cdot\nabla\phi\right](Bx)dx,
\end{align}
for any $\phi\in C_{c}^1(\mathbb{R}^{n})$, and $\widetilde{\phi}(x)=\phi(Bx)$. In other words, in the weak sense, there holds $\Delta_p \widetilde u = (\Delta_p^H u)(Bx)$.
\end{thm}

\begin{rem}
The aforementioned result indeed means that under the assumption (\ref{hm}), for any $p>1$, the anisotropic equation can be transformed into an isotropic equation via a linear transformation. Specifically, if $u$ is a weak solution of equation (\ref{main eq1}), then $\widetilde{u}$ serves as a weak solution to the equation
$$-\Delta_{p}\widetilde{u} = g(\widetilde{u}) \quad \mbox{in }\; \mathbb{R}^{n}.$$
\end{rem}

As direct application of Theorem \ref{thm1}, we show quickly the following relationships between the classical Kelvin transform and the anisotropic Kelvin transform, under the assumption (\ref{hm}).

\begin{thm}\label{thm2}
Suppose that $H$ satisfies (\ref{hm}) and $u \in H^1_{loc}(\R^n\backslash\{0\})$. Let $\widehat{u}:=\frac{u\circ T_{H}}{H^{n-2}}$ and $\eta=Bx$ with $B = \sqrt{M}$, then
$$\Delta^{H^*}\widehat{u}(x) = \Delta_{\eta} [\widehat{u}(B^{-1}\eta)] =\frac{(\Delta^{H}u)\circ T_{H}}{H^{n+2}}(x) \quad \mbox{in }\; \R^n\backslash\{0\}.$$
Here $\Delta_\eta$ means the Laplacian w.r.t.~the new variable $\eta$. Similar result holds for $\Delta_n^H$ in $\R^n$, that is, for any $u \in W^{1, n}_{loc}(\R^n\backslash\{0\})$, set $u^{*}:=u\circ T_{H}$, there holds 
$$\Delta_n^{H^*}u^*(x) = (\Delta_{n})_{\eta}[u^{*}(B^{-1}\eta)]=\frac{(\Delta_{n}^{H}u)\circ T_{H}}{H^{2n}}(x) \quad \mbox{in }\; \R^n\backslash\{0\}.$$
\end{thm}

\begin{rem}
Readily, we can deduce Propositions A and B by Theorems \ref{thm2}. Moreover, we note that several assumptions in Propositions A and B, like local boundedness of $u$, $f \in L^2_{loc}$ or $g \in L^\frac{n}{n-1}_{loc}$ are not relevant.
\end{rem}

As previously elucidated, under the assumption $(H_M)$, the operator $\Delta_p^H$ acts like the classical $p$-Laplacian up to a linear transform, in particular the Finsler Laplacian is indeed a linear elliptic operator with constant coefficients. Owing to the criterion established by Cozzi, Farina, and Valdinoci (Proposition D), it is now apparent that the average formulae proposed by Ferone and Kawohl (Proposition C) originate from the mean value property of harmonic functions. Along similar lines, numerous other findings pertaining to equation \eqref{main eq1} under the assumptions (\ref{hm}) or (\ref{hm1}) can be expediently derived from analogous considerations involving the isotropic $p$-Laplacian. For instance, the principal results presented in \cite{FL2021,FL2022} can be achieved by applying previous research conducted in \cite{DF2009, W2012}.

\smallskip
Consider the classification problem with the weighted Finsler Liouville equation
\begin{equation}
\label{WFL}
-\Delta^{H}u=[H^{*}(x)]^{\alpha}e^{u} \;\; \mbox{ in }\;\mathbb{R}^2,\quad \int_{\mathbb{R}^2}[H^{*}(x)]^{\alpha}e^{u} dx<+\infty.
\end{equation}
When $\alpha=0$, this question was initially posed by Wang and Xia \cite{WX2012}. More recently, Ciraolo and Li \cite{CL2023} have demonstrated that, up to translation and scaling, any weak solution to \eqref{WFL} takes the form
$$u(x)=-2\ln\left[1+\frac{H^{*}(x)^2}{8}\right].$$
In the context of anisotropy, we refer to functions that solely depend on $H^{*}$ as radial functions.

\smallskip
For the general case where $\alpha>-2$, the solutions of equation \eqref{WFL} become significantly more intricate. Indeed, for the isotropic case, $H(\xi)=|\xi|$, with $\Delta^H=\Delta$, Prajapat and Tarantello \cite{PT2001} employed methods such as the Kelvin transform, moving plane technique, and complex analysis to fully classify the solutions to \eqref{WFL}. They found that for $\alpha\notin\mathbb{N}$ or $\alpha=0$, any solution is radially symmetric about the origin, whereas for $\alpha\in\mathbb{N}\setminus\{0\}$, non-radially symmetric solutions exist. These findings suggest that a similar classification result may hold when $H$ satisfies the $(H_M)$ condition. However, an intriguing question arises: can we generalize this classification result to the general anisotropic operator $\Delta^H$ without the assumption \eqref{hm}?

\section{Proof of main results}

In this section, we will prove our main claims through basic algebraic operations.

\subsection{Proof of Theorem \ref{thm1}.}
Assume that $M = (m_{ij})_{1\leq i, j \leq n}$ is a symmetric and positive definite matrix and $H$ is given by \eqref{hm}. For any $u \in C^1(\R^n)$, $\phi\in C_{c}^1(\mathbb{R}^{n})$, straightforward calculation yields
\begin{align*}
H^{p-1}(\nabla u)H_{\xi}(\nabla u)\cdot\nabla\phi &= H^{p-2}(\nabla u)\sum_{i,j=1}^{n}m_{ij}u_{x_{j}}\phi_{x_{i}}
\end{align*}
where $u_{x_{i}}=\frac{\partial u}{\partial x_{i}}$. Let $\sqrt{M} = B = (b_{ij})_{1\leq i, j \leq n}$ and $\widetilde{u}(x):=u(Bx)$, then
\begin{align}
\widetilde{u}_{x_{i}}(x) =\sum_{j=1}^{n}b_{ij}u_{y_{j}}(Bx),\nonumber
\end{align}
and
\begin{align*}
|\nabla \widetilde{u}|^2(x) = \sum_{i=1}^{n}\widetilde{u}_{x_{i}}^{2}(x) = \sum_{i,j,k=1}^{n}b_{ij}b_{ik}u_{y_{j}}(Bx)u_{y_{k}}(Bx) & =\sum_{j,k=1}^{n}m_{jk}u_{y_{j}}(Bx)u_{y_{k}}(Bx)\\
& =\left[H(\nabla u)\right]^2(Bx).
\end{align*}
Similarly, for $\widetilde{\phi}(x):=\phi(Bx)$, there holds
$$\widetilde{\phi}_{x_{i}}(x) =\sum_{j=1}^{n}b_{ij}\phi_{y_{j}}(Bx),$$
and
$$\nabla\widetilde{u}\cdot \nabla\widetilde{\phi} (x) =\sum_{i,j,k=1}^{n}b_{ij}b_{ik}\phi_{y_{j}}(Bx)u_{y_{k}}(Bx)=\sum_{i,j=1}^{n}m_{ij}u_{y_{i}}(Bx)\phi_{y_{j}}(Bx).$$
It follows that
\begin{align*}
|\nabla \widetilde{u}|^{p-2}\nabla\widetilde{u}\cdot\nabla\widetilde{\phi} &= \left[H(\nabla u)\right]^{p-2}(Bx)\sum_{i,j=1}^{n}m_{ij}u_{y_{i}}(Bx)\phi_{y_{j}}(Bx)\\
&= \left[\left[H(\nabla u)\right]^{p-1}H_{\xi}(\nabla u)\cdot\nabla\phi\right](Bx).
\end{align*}
Hence for $u \in C^1(\R^n)$, we get \eqref{equa1}, or equivalently $\Delta_p \widetilde u = (\Delta_p^H u)(Bx)$ in the sense of Definition \ref{def1}. Applying standard density argument, the result remains valid for all $u \in W^{1, p}_{loc}(\R^n)$. $\hfill\square$

\subsection{Proof of Theorem \ref{thm2}.}
Let $H$ be given by \eqref{hm}. Denote $B = \sqrt{M}$ and $L_B(x) = Bx$. Since $H^{2}(x) = x^tMx = |Bx|^2$, we see that
$$T_H(x) = \frac{\nabla H(x)}{H(x)} = \frac{B^2x}{|Bx|^2} = L_B\circ \mathcal{I}\circ L_B(x),$$
where $\mathcal{I}(x) = \frac{x}{|x|^2}$ is the spherical inversion in Euclidean spaces.

Moreover, let $\eta = L_B(x) = Bx$, as $M^{-1} = (B^{-1})^2$ and $H^*$ satisfies the $(H_M)$ assumption with $M^{-1}$, applying Theorem \ref{thm1} with $p = 2$, we have
\begin{align*}
\Delta^{H^{*}} \widehat{u}(x) = \Delta_\eta \left[ \widehat{u}(B^{-1}\eta)\right].
\end{align*}
As $H(x) = |Bx| = |\eta|$ with the new variable $\eta$, there holds
$$\widehat{u}(B^{-1}\eta) = \frac{u\circ T_H(B^{-1}\eta)}{|\eta|^{n-2}} = \frac{u\circ L_B\circ \mathcal{I}(\eta)}{|\eta|^{n-2}}=  \frac{(u\circ L_B)\circ \mathcal{I}(\eta)}{|\eta|^{n-2}} = (u\circ L_B)_{\mathcal K},$$
where $w_{\mathcal K}$ means the classical Kelvin transform in Euclidean spaces for the function $w$. By the standard relationship between $\Delta w_{\mathcal K}$ and $\Delta w$ in $\R^n\backslash\{0\}$, we get
\begin{align*}
\Delta^{H^{*}} \widehat{u}(x) = \Delta_\eta \left[\widehat{u}(B^{-1}\eta)\right] = \Delta_\eta (u\circ L_B)_{\mathcal K} &= |\eta|^{-2-n} \left[\Delta_\eta (u\circ L_B)\right]\circ \mathcal{I}\\
& = |\eta|^{-2-n} (\Delta^H u)\circ L_B\circ \mathcal{I}(\eta)\\
& = |Bx|^{-2-n} (\Delta^H u)\circ T_H(x)\\
& = \frac{(\Delta^H u)\circ T_H(x)}{H^{n+2}(x)}.
\end{align*}
Here we applied Theorem \ref{thm1} to get the second line.

\medskip
Moreover, recalling that $H_{\xi}(y) = \frac{My}{H(y)}$ for any $y \ne 0$, it follows from Theorem \ref{thm1} that
$$\Delta_{n}^{H^{*}}u(x)=(\Delta_{n})_{\eta}[u(B^{-1}\eta)].$$
For the same reason, let $x = B\zeta$, we have
\begin{align*}
\Delta_{n}^{H}u(x) =(\Delta_{n})_{\zeta}[u\circ L_B(\zeta)].
\end{align*}
Therefore
\begin{align*}
\Delta_{n}^{H^{*}}u^*(x) = (\Delta_{n})_{\eta}[u^{*}(B^{-1}\eta)] &=(\Delta_{n})_{\eta}\left[(u\circ L_B)\circ\mathcal{I}\right] \\
& =|\eta|^{-2n}[(\Delta_{n})_{\eta}(u\circ L_B)]\circ\mathcal{I}\\
&=|\eta|^{-2n}(\Delta_{n}^{H}u)\circ T_{H}(x)=\frac{(\Delta_{n}^{H}u)\circ T_{H}}{H^{2n}}(x).
\end{align*}
The second equality above comes from the Kelvin transform for $\Delta_n$ in $\R^n$. $\hfill\square$

\smallskip


\smallskip


\section*{Acknowledgements}
The authors would like to thank the anonymous referees for their careful reading. Y. Li is partially supported by China Postdoctoral Science Foundation (No. 2022M721164). Both authors are supported in part by Science and Technology Commission of Shanghai Municipality (No. 22DZ2229014).

\end{document}